\documentclass[10pt,prd,longbibliography,superscriptaddress]{revtex4-2}

\usepackage[T1]{fontenc}
\usepackage[utf8]{inputenc}
\usepackage[english]{babel}

\usepackage{amsfonts,amsbsy,amssymb,amsmath,graphicx,float} 
\usepackage{color}

\graphicspath{{figures/}}

\usepackage[linktoc=all,colorlinks=true,citecolor=blue,linkcolor=blue]{hyperref}

\begin{document}

\title{The Influence of a Parameter that Controls the Asymmetry of a Potential Energy Surface with an Entrance Channel and Two Potential Wells}

\author{Makrina Agaoglou}
\email{makrina.agaoglou@bristol.ac.uk}
\affiliation{School of Mathematics, University of Bristol, \\ Fry Building, Woodland Road, Bristol, BS8 1UG, United Kingdom.}
\affiliation{Instituto de Ciencias Matematicas, C/ Nicolas Cabrera 13-15, Campus Cantoblanco, 28049 Madrid, Spain.}

\author{Matthaios Katsanikas}
\email{matthaios.katsanikas@bristol.ac.uk}
\affiliation{School of Mathematics, University of Bristol, \\ Fry Building, Woodland Road, Bristol, BS8 1UG, United Kingdom.}
\affiliation{Research Center for Astronomy and Applied Mathematics, Academy of Athens, Soranou Efesiou 4, Athens, GR-11527, Greece.}

\author{Stephen Wiggins}
\email{s.wiggins@bristol.ac.uk}
\affiliation{School of Mathematics, University of Bristol, \\ Fry Building, Woodland Road, Bristol, BS8 1UG, United Kingdom.}

\begin{abstract}
In this paper we study an asymmetric valley-ridge inflection point (VRI) potential, whose energy surface (PES) features  two sequential index-1 saddles (the upper and the lower), with one saddle having higher energy than the other and two potential wells separated by the lower index-1 saddle. We show how the depth and the flatness of our potential changes as we modify the parameter that controls the asymmetry as well as how the branching ratio (ratio of the trajectories that enter each well) is changing as we modify the same parameter and its correlation with the area of the lobes as they have formulated by the stable and unstable manifolds that have been extracted from the gradient of the LD scalar fields.  


\end{abstract}

\maketitle


\section{Introduction}
\label{sec:intro}
In this paper we consider a two degree-of-freedom Hamiltonian system having a valley ridge inflection point (VRI) potential energy surface (PES).  VRI potential energy surfaces have four critical points: a high energy saddle and a lower energy saddle separating two wells. In between the two saddle points is a valley ridge inflection point that is the point where the potential energy surface geometry changes from a valley to a ridge. The region between the two saddles forms a reaction channel and the dynamical issue of interest is how trajectories cross the high energy saddle, evolve towards the lower energy saddle, and select a particular well to enter. This selectivity is controlled by the branching ratio. The goal of this paper is to analyze the effect of asymmetry of the two potential wells on the branching ratio (the study of the effects of the asymmetry in a potential energy surface is important in chemical reaction dynamics, dynamical astronomy and Hamiltonian fluid dynamics). For this purpose we use a simplified version of the VRI PES that was first discussed in \cite{collins2013}) and that allows the symmetry to be controlled with a single parameter. The properties of the PES are further quantified using the notions of depth and flatness of the PES introduced in \cite{naikflat}.

\begin{figure}[htbp]
	\begin{center}
		\includegraphics[scale=0.60]{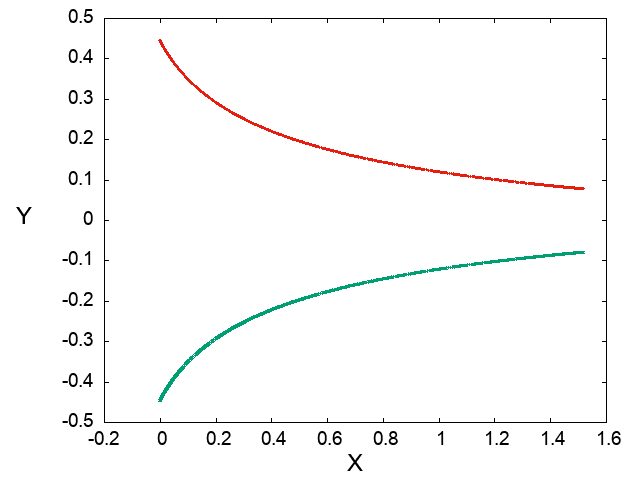} 
	\end{center}
	\caption{The top (with red) and bottom (with green) unstable periodic orbits in the configuration space for  $E=0.1,c=0$. }
	\label{per-o}
	\end{figure}

In previous papers 
\cite{Agaoglou2020,Katsanikas2020,crossley2021poincare} we have studied the phase 
space mechanism for selectivity, as quantified by a branching ratio, in the symmetric case of this PES. For this symmetric PES the  fraction of trajectories that enter each well given an initial number of trajectories is $1 : 1$. Moreover in the aforementioned papers, we present how  the selectivity is a consequence of the heteroclinic and homoclinic  connections established between the invariant manifolds of the  families of unstable periodic orbits (UPOs) (see the top and bottom unstable periodic orbits, that are presented in the system, in the configuration space in Fig. \ref{per-o})  present in the system. Furthermore, in the paper \cite{douglas2021} we studied the time evolution of the trajectories after the selectivity. A detailed study of the bifurcations of the periodic orbit dividing surfaces in this system is presented in \cite{katsanikas2021bifurcation}. In the paper \cite{Victor2020} we consider the same symmetric PES but now we  add to the Hamiltonian model a time-periodic forcing term. This  forcing term depends on an amplitude, frequency, and phase and our goal was to analyse how the branching ratio depends on these three parameters. 

In this paper  we study the asymmetric case of a potential  with one saddle having higher energy than the other and two potential wells separated by the lower index-1 saddle. We explore the effects of the asymmetry parameter in the depth and the flatness of our potential. Moreover we study the phase space structure of the asymmetric case. In particular, we use the method of Lagrangian descriptors (LD)  to compute and visualize the invariant manifolds of the unstable periodic orbits that are responsible for the transport of the trajectories from the region of the upper index-1 saddle to the region of the wells. We did this in order to understand the correlation between the phase space structure and the branching ratio in the asymmetric case of our potential. 

The method of Lagrangian Descriptors has been introduced more than a decade ago and it is a trajectory-based scalar technique. This technique has the ability to reveal the geometrical template of phase space structures and has been used in several applications in different scientific areas i.e. Geophysics \cite{mendoza2010,curbelo2019}, Chemistry \cite{Agaoglou2020,Katsanikas2020,gonzalez2020}. 

This paper is organized as follows. In Section \ref{Model} we present the two degrees of freedom (DoF) Hamiltonian model that we used for our analysis. In Section \ref{LDs} we present briefly the method of Lagrangian Descriptors. Initially in section \ref{Results} we investigate  the depth and the flatness of our potential and then we provide a detailed description of our set up and show how the phase space structure of our system is related to the branching ratio of the trajectories that choose one well over the other. Finally in Section \ref{Conclusions} we summarize our results.

\section{Model}\label{Model}

\begin{figure}[htbp]
	\begin{center}
    	A)\includegraphics[scale=0.20]{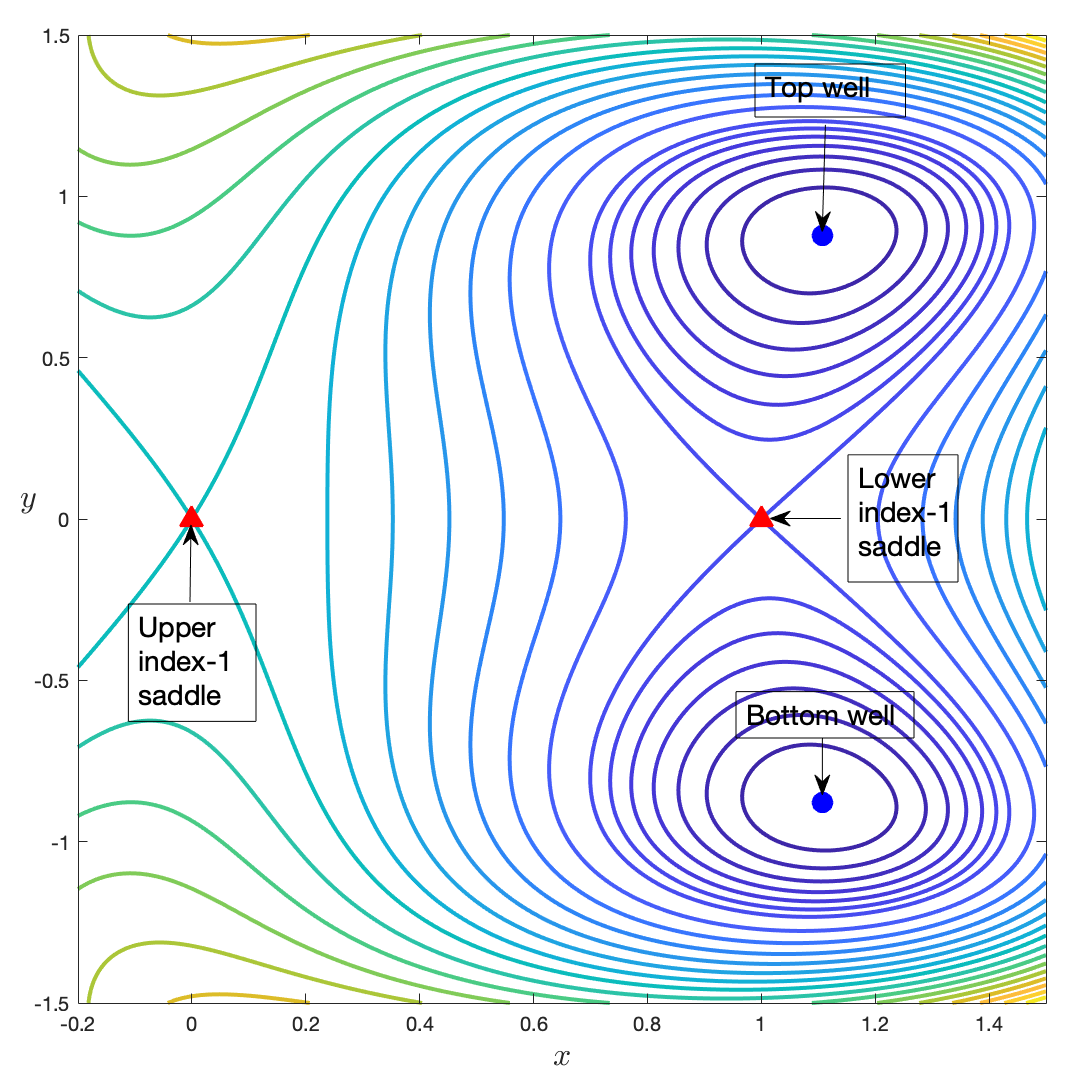} 
		B)\includegraphics[scale=0.20]{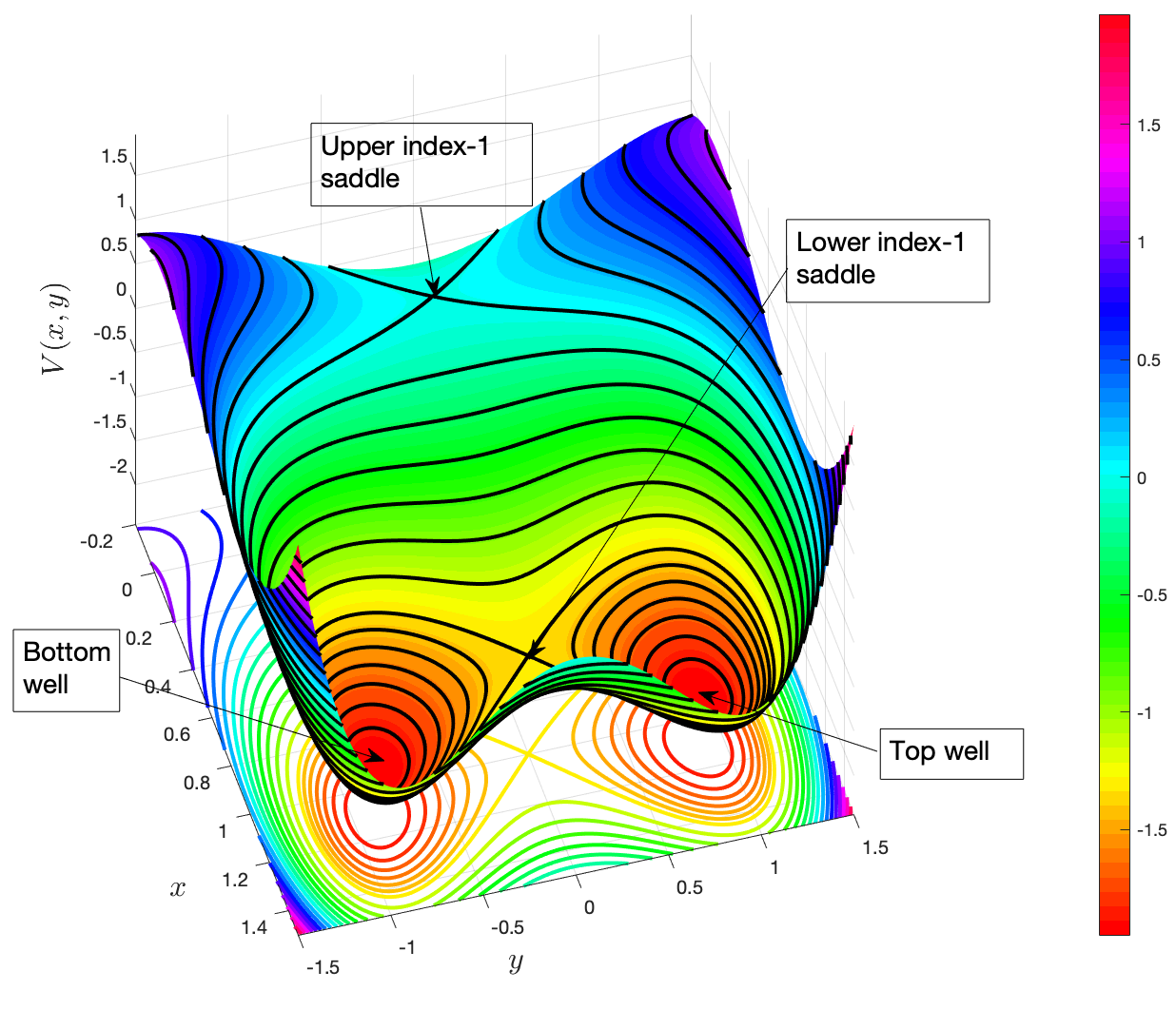} \\
		C)\includegraphics[scale=0.20]{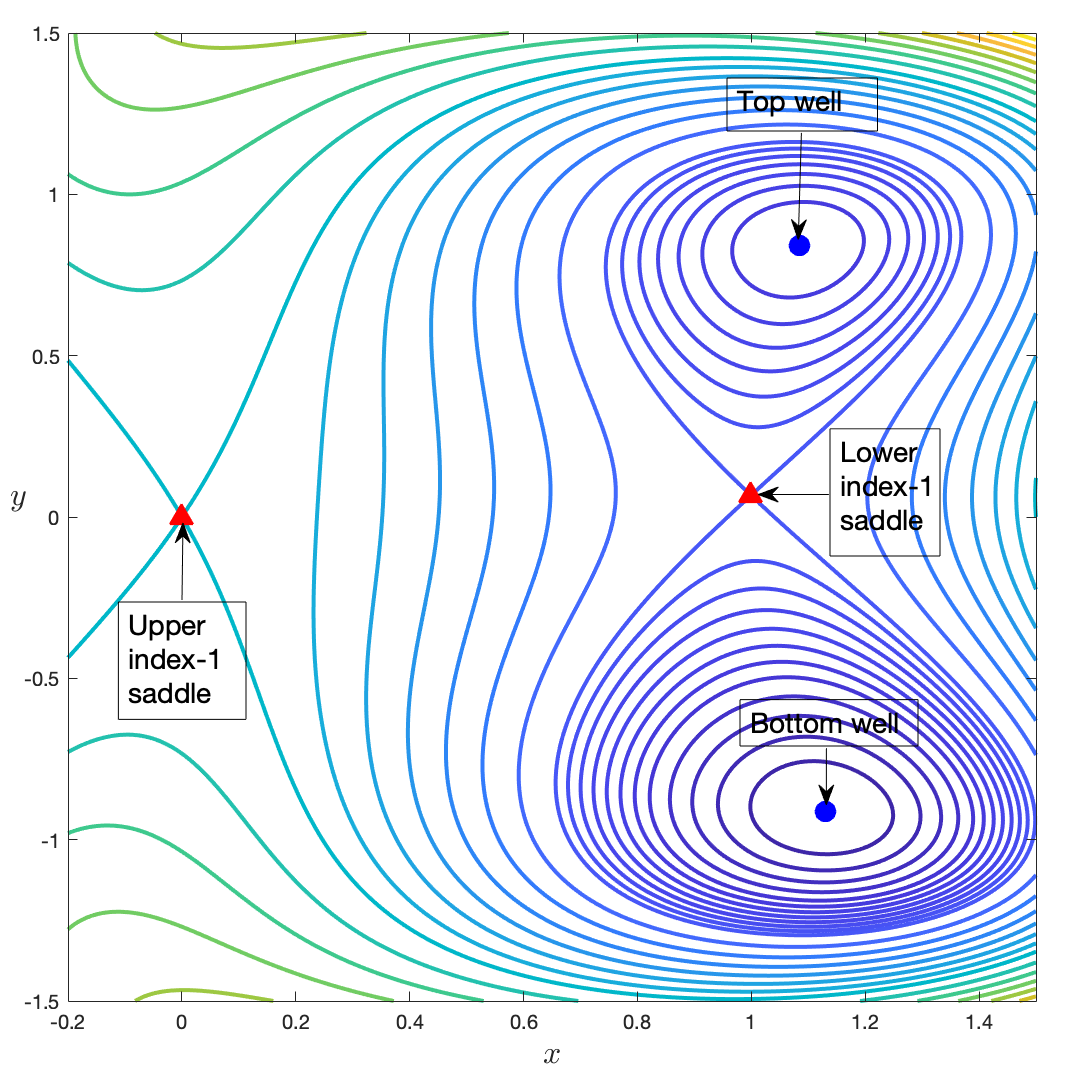} 
		D)\includegraphics[scale=0.20]{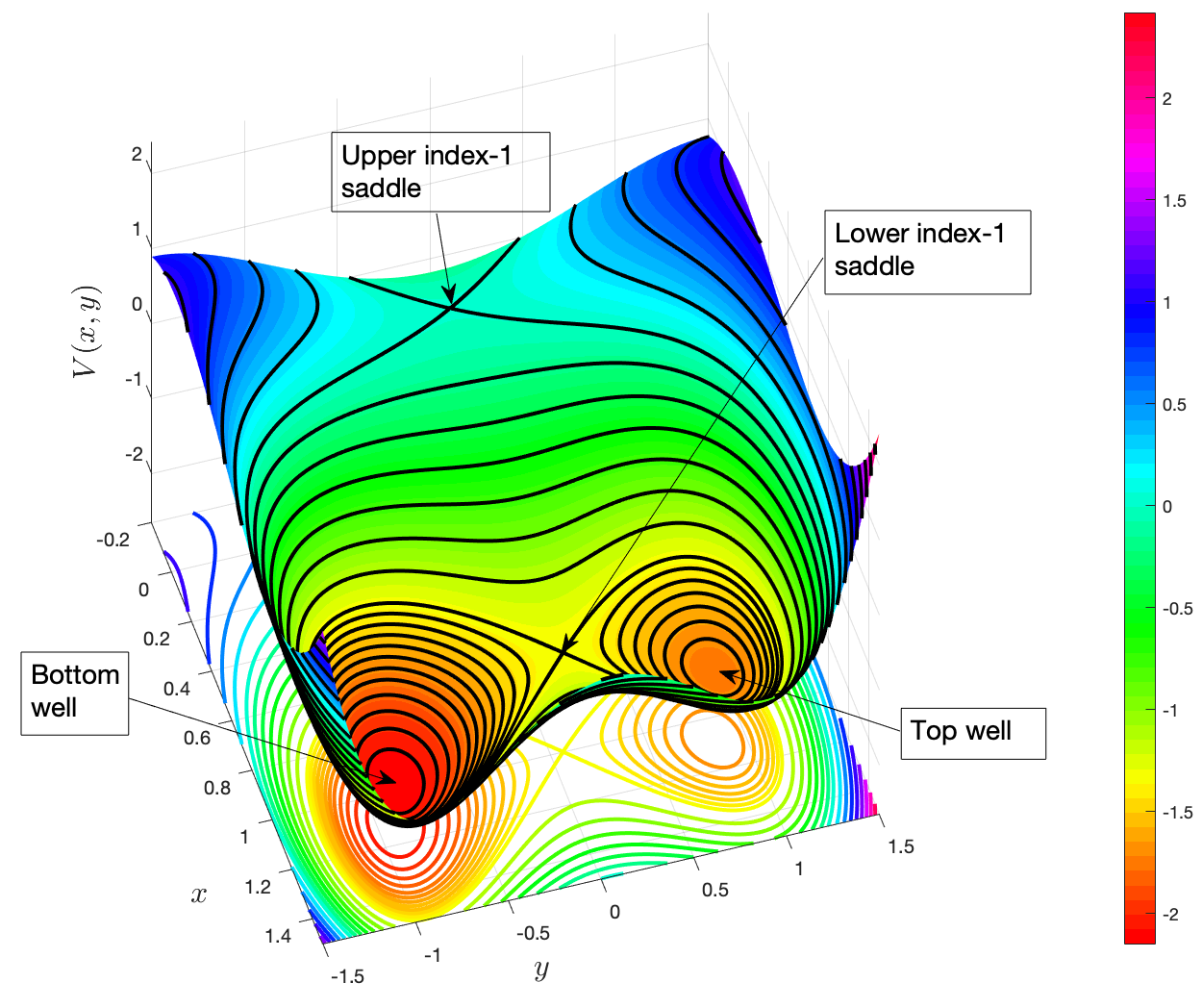} \\
	    E)\includegraphics[scale=0.215]{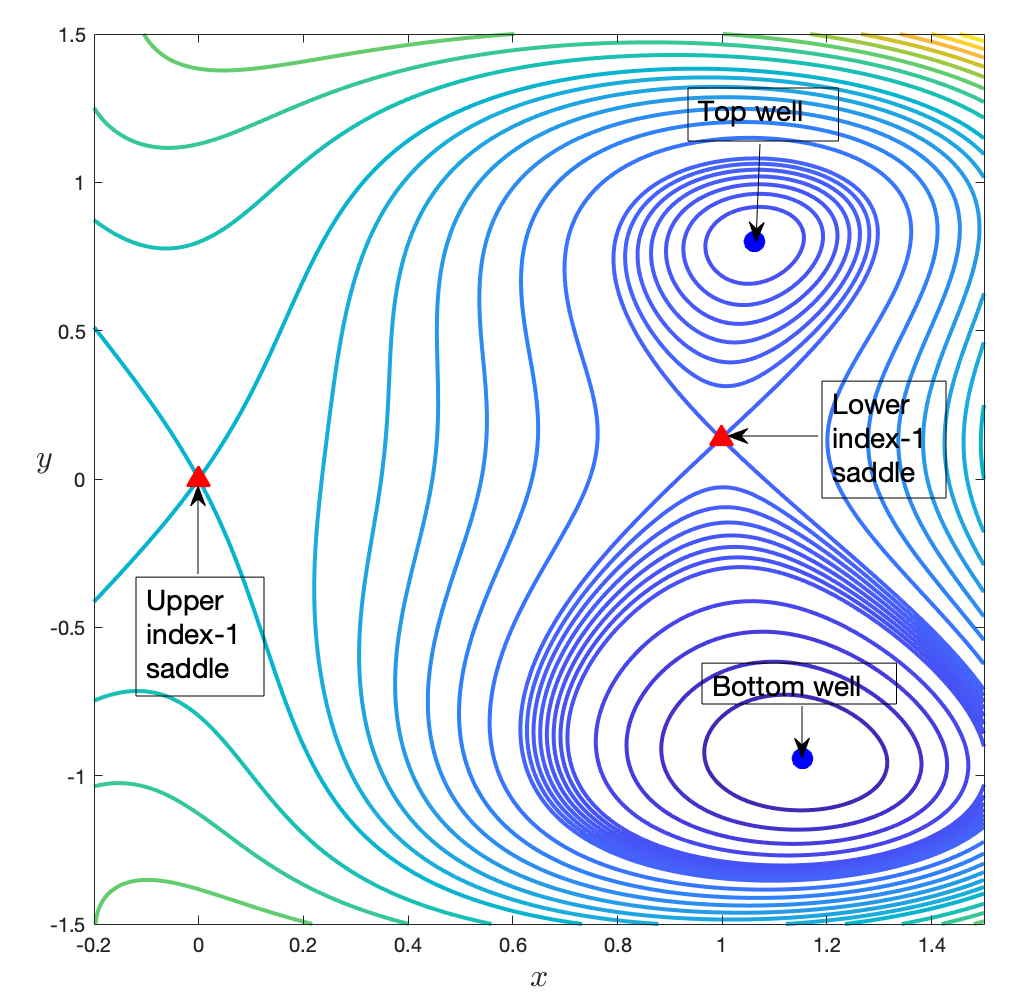} 
		F)\includegraphics[scale=0.205]{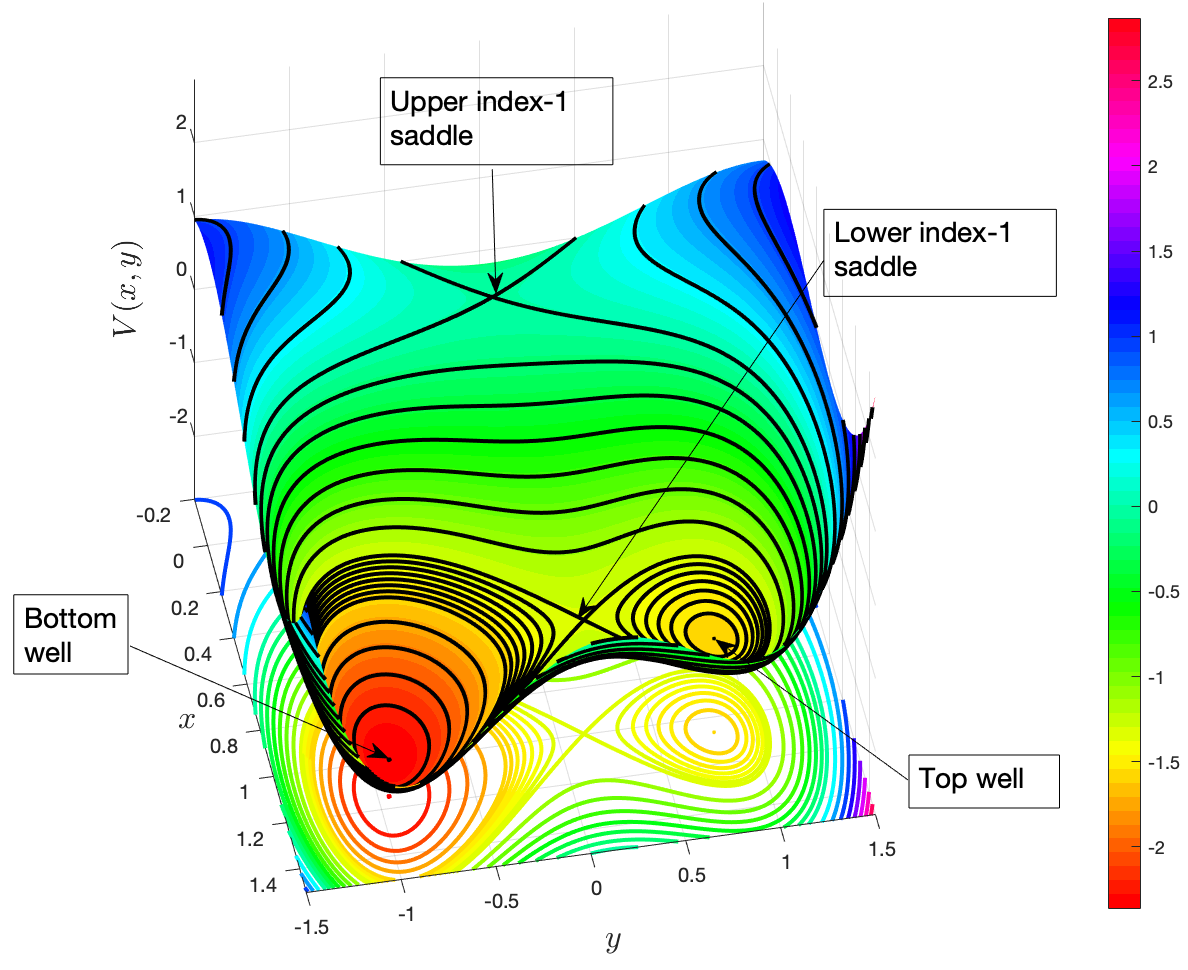} 
	\end{center}
	\caption{PES for different values of the parameter $c$: A-B) $c = 0$ (symmetric case), C-D) $c = 0.2$ and E-F) $c = 0.4$.}
	\label{fig:PES}
\end{figure}

In this paper we use a simplified version of the PES from the one discussed in \cite{collins2013}. Our Hamiltonian model is the sum of kinetic plus potential energy of the form:

\begin{equation}
H(x,y,p_x,p_y) = \dfrac{p_x^2}{2 m_x} + \dfrac{p_y^2}{2 m_y} + V(x,y) \;,
\label{hamiltonian}
\end{equation}

\noindent
where we consider that the mass in each DoF is $m_x = m_y = 1$, and the PES is given by:

\begin{equation}
V(x,y) = \dfrac{8}{3}x^3 - 4x^2 + \dfrac{1}{2} y^2 + x\left(y^4 - 2 y^2\right)+ cxy 
\label{pes_modelVRI}
\end{equation}  

\noindent
where the parameter $c$ controls the asymmetry of the wells, as we will  see below, with respect to the x axis. The PES has two wells separated by an index-1 saddle. We refer to this saddle as the lower index-1 saddle. Furthermore the PES has an entrance/exit channel determined by an index-1 saddle located at the origin. We refer to this saddle as the upper saddle. The upper saddle has the highest energy of all the critical points on the PES.  In particular, we illustrate the PES in Fig. \ref{fig:PES} for three different values of the parameter $c$: $c = 0$ (the symmetric case), $c = 0.2$ and $c = 0.4$ and we also illustrate the location and the energies of all the critical points in Tables \ref{tab:tab3}, \ref{tab:tab2},\ref{tab:tab1}, respectively. It is easily observed that in the symmetric case the energy of both of the wells is equal but as we increase the value of the asymmetry parameter, $c$, the top well becomes more flat than the bottom well and therefore its energy is higher than the energy of the bottom well.

Hamilton's equations of motion for our model as given in Eq. (\ref{pes_modelVRI}) are as follows:

\begin{equation}
\begin{cases}
\dot{x} = \dfrac{\partial H}{\partial p_x} = p_x \\[.5cm]
\dot{y} = \dfrac{\partial H}{\partial p_y} = p_y  \\[.5cm]
\dot{p}_x = -\dfrac{\partial H}{\partial x} = 8 x \left(1 - x\right) + y^2\left(2 - y^2 \right) - cy \\[.5cm]
\dot{p}_y = -\dfrac{\partial H}{\partial y} = y \left[4 x \left(1 - y^2\right) - 1\right] - cx
\end{cases}
\;.
\label{ham_eqs}
\end{equation}

\begin{table}[htbp]	
	\begin{tabular}{| l | c | c | c | c |}
		\hline
		Critical point \hspace{1cm} & \hspace{0.6cm} x \hspace{0.6cm} & \hspace{.6cm} y \hspace{.6cm} & \hspace{.2cm} \text{Potential Energy} $(V)$ \hspace{.2cm} & \hspace{.6cm} \text{Stability} \hspace{.6cm} \\
		\hline\hline
		index-1 saddle (Upper) \hspace{.5cm} & 0 & 0 & 0 & saddle $\times$ center \\
		\hline
		index-1 saddle (Lower) \hspace{.5cm} & 1 & 0 & -4/3 & saddle $\times$ center \\
		\hline
		Potential Well (Top) \hspace{.5cm} & 1.1071 & 0.8799  & -1.9477 & center \\
		\hline
		Potential Well (Bottom) \hspace{.5cm} & 1.1071 & -0.8799 &  -1.9477 & center  \\
		\hline
		\end{tabular} 
		\caption{Location of the critical points of the potential energy surface for $c=0$.} 
	\label{tab:tab3} 
\end{table}

\begin{table}[htbp]	
	\begin{tabular}{| l | c | c | c | c |}
		\hline
		Critical point \hspace{1cm} & \hspace{0.6cm} x \hspace{0.6cm} & \hspace{.6cm} y \hspace{.6cm} & \hspace{.2cm} \text{Potential Energy} $(V)$ \hspace{.2cm} & \hspace{.6cm} \text{Stability} \hspace{.6cm} \\
		\hline\hline
		index-1 saddle (Upper) \hspace{.5cm} & 0 & 0 & 0 & saddle $\times$ center \\
		\hline
		index-1 saddle (Lower) \hspace{.5cm} & 0.9994 & 0.0671 & -1.3266 & saddle $\times$ center \\
		\hline
		Potential Well (Top) \hspace{.5cm} & 1.0845 & 0.8427 & -1.7587 & center  \\
		\hline
		Potential Well (Bottom) \hspace{.5cm} & 1.1303 & -0.9130 & -2.1481 & center \\

		\hline
		\end{tabular} 
		\caption{Location of the critical points of the potential energy surface for $c=0.2$.} 
	\label{tab:tab2} 
\end{table}

\begin{table}[htbp]	
	\begin{tabular}{| l | c | c | c | c |}
		\hline
		Critical point \hspace{1cm} & \hspace{0.6cm} x \hspace{0.6cm} & \hspace{.6cm} y \hspace{.6cm} & \hspace{.2cm} \text{Potential Energy} $(V)$ \hspace{.2cm} & \hspace{.6cm} \text{Stability} \hspace{.6cm} \\
		\hline\hline
		index-1 saddle (Upper) \hspace{.5cm} & 0 & 0 & 0 & saddle $\times$ center \\
		\hline
		index-1 saddle (Lower) \hspace{.5cm} & 0.9978 & 0.1368 & -1.3064 & saddle $\times$ center \\
		\hline
		Potential Well (Top) \hspace{.5cm} & 1.0624 & 0.7998 & -1.5819 & center  \\
		\hline
		Potential Well (Bottom) \hspace{.5cm} & 1.1541 & -0.9431 & -2.3592 & center \\
		\hline
		\end{tabular} 
		\caption{Location of the critical points of the potential energy surface for $c=0.4$.} 
	\label{tab:tab1} 
\end{table}

\section{Lagrangian Descriptors}\label{LDs}

The method of Lagrangian Descriptors (LDs) is a 
trajectory-based scalar diagnostic. It was first introduced to analyze Lagrangian transport and mixing in Geophysical flows 
\cite{madrid2009,mendoza2010}. This methods explores the geometrical template of phase space structures. There are several definitions of  LDs. The first one introduced was based on the computation of the arclength of the trajectories of initial conditions as they evolve forward and backward in time \cite{mendoza2010,mancho2013lagrangian} for a fixed integration time.

We consider the following dynamical system with time dependence:

\begin{equation}
\dfrac{d\mathbf{x}}{dt} = \mathbf{v}(\mathbf{x},t) \;,\quad \mathbf{x} \in \mathbb{R}^{n} \;,\; t \in \mathbb{R} \;
\label{eq:gtp_dynSys}
\end{equation}

\noindent
where $\mathbf{v}(\mathbf{x},t) \in C^{r} (r \geq 1)$ in $\mathbf{x}$ and it is continuous in time.

We will introduce an alternative definition of LDs, the $p-$norm definition, which was first presented in \cite{lopesino2017} and it is the one that we use in this paper. Given an initial condition $x_0$ at time $t_0$, take a fixed integration time $\tau>0$ and $p \in (0, 1]$. The definition of the LDs is the following:

\begin{equation}
M_p(\mathbf{x}_{0},t_0,\tau) = \sum_{k=1}^{n} \bigg[ \int^{t_0+\tau}_{t_0-\tau}  |v_{k}(\mathbf{x}(t;\mathbf{x}_0),t)|^p \; dt \bigg] = M_p^{(b)}(\mathbf{x}_{0},t_0,\tau)+ M_p^{(f)}(\mathbf{x}_{0},t_0,\tau) \;, 
\label{eq:Mp_function}
\end{equation}

\noindent
where $M_p^{(b)}$ and $M_p^{(f)}$ are its backward and forward integration parts:

\begin{equation}
\begin{split}
M_p^{(b)}(\mathbf{x}_{0},t_0,\tau) &= \sum_{k=1}^{n} \bigg[ \int^{t_0}_{t_0-\tau}  |v_{k}(\mathbf{x}(t;\mathbf{x}_0),t)|^p \; dt \bigg] \;, \\[.2cm]
M_p^{(f)}(\mathbf{x}_{0},t_0,\tau) &= \sum_{k=1}^{n} \bigg[ \int^{t_0+\tau}_{t_0} |v_{k}(\mathbf{x}(t;\mathbf{x}_0),t)|^p \; dt \bigg].
\end{split}
\end{equation}

In our computations we consider $p=1/2$. The forward integration reveals the stable manifolds of our dynamical system whereas the backward integration reveals the unstable manifolds of our dynamical system and therefore, by combining both reveals all the invariant manifolds. Moreover, in this paper all the analysis in the phase space structures has been carried out in the Poincare section:

\begin{equation}
    \Sigma (H_{0}) = \{(x,y,p_{x},p_{y})\in \mathbb{R}^{4}| x = 0.05, p_{x}(x,y,p_{y};H_{0})>0\}
\end{equation}

\section{Results}\label{Results}

In this section, we describe our results concerning the effects of the asymmetry parameter on the depth and flatness of our potential energy surface (subsection A). Furthermore, we describe how the the area of the lobes that are formed between  the invariant manifolds which are responsible for the transport of trajectories from the region of the upper index-1 saddle to the wells, is related  to the branching ratio of the trajectories that choose one well over the other (subsection B).

\subsection{Depth and flatness of the Valley Ridge Inflection Point Potential Energy Surface}

\begin{figure}[htbp]
	\begin{center}
		A)\includegraphics[scale=0.20]{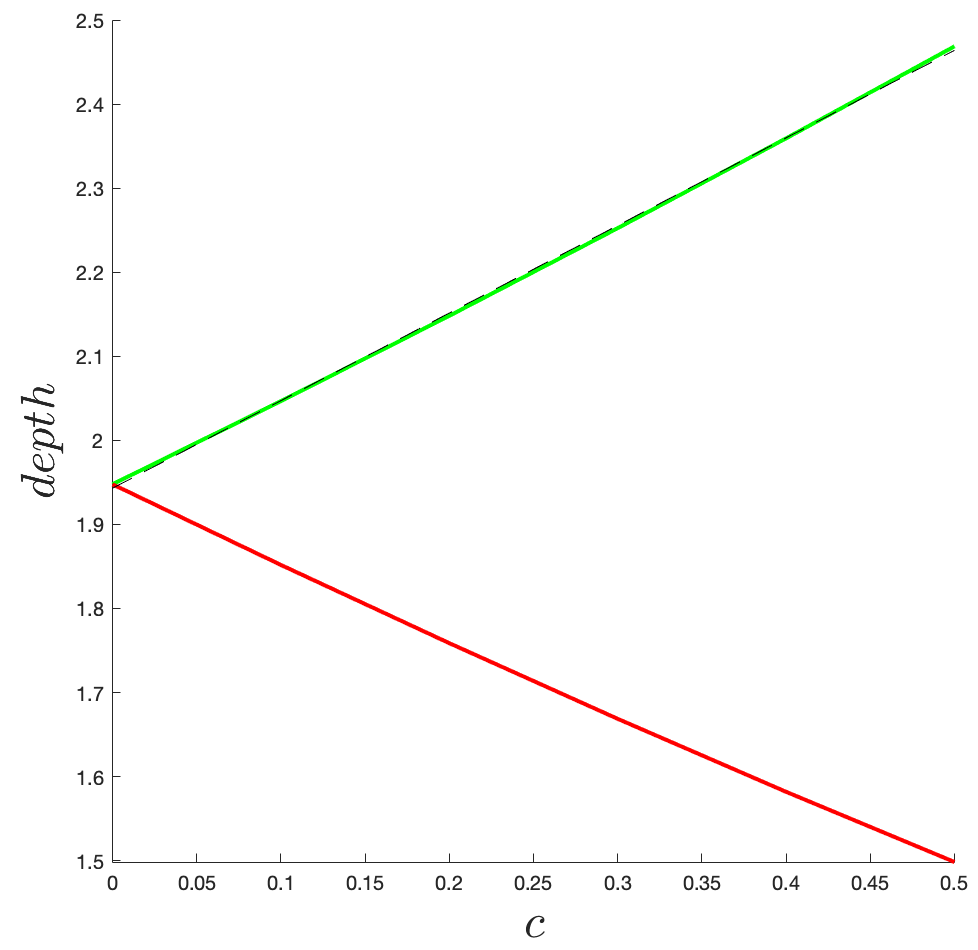} 
		B)\includegraphics[scale=0.20]{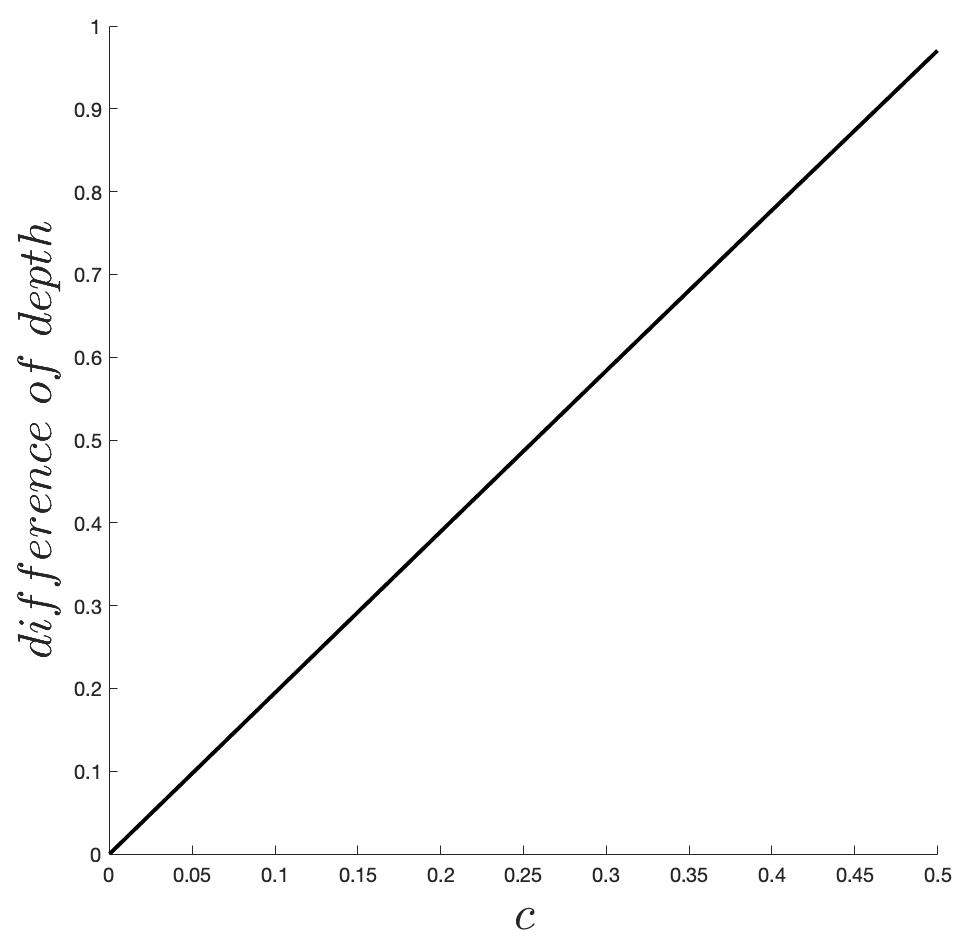}\\
		C)\includegraphics[scale=0.20]{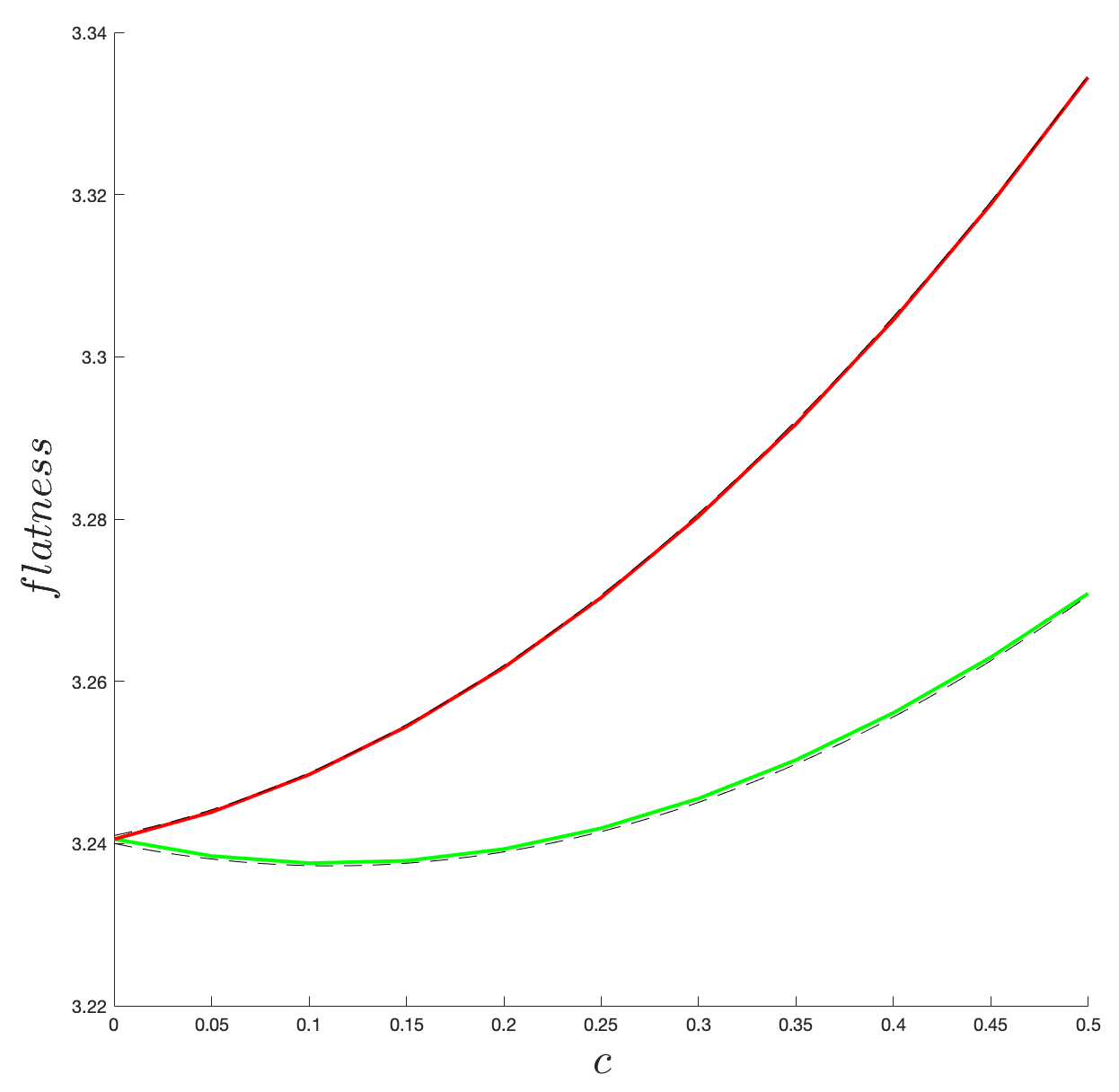}
	\end{center}
	\caption{A) The depth of the PES for our system when the center equilibrium point is located in the bottom (top) well is indicated with a green (red) color. Superimposed with the linear approximation $y = d_{1}c + d_{2}$, where $d_{1} = 1.042, d_{2} =1.943$ (black dashed line). B) The difference between the depth of the two wells. C) The flatness for our system. In red color we present the flatness in the bounded domain $[0.6,1.5]\times [-1.5,0]$ and in green color we present the flatness in the bounded domain $[0.6,1.5]\times [0,1.5]$ superimposed with the quadratic approximations $h = f_{1}c^{2} + f_{2}c + f_{3}$ and $k = f_{4}c^{2} - f_{5}c + f_{6}$ respectively (black dashed lines), where $f_{1} =0.2749, f_{2} = 0.04992, f_{3} = 3.241$ and $f_{4} =0.2207, f_{5} =0.04926 , f_{6} = 3.24$. The asymmetry parameter varies between 0 and 0.5.}
	\label{depth_flat}
\end{figure}

In this subsection we study the depth and flatness of the VRI PES of our system. In the paper \cite{naikflat} they show how the geometry of a potential energy surface with one well may change as a function of the depth and the flatness and how these changes affect the reaction dynamics. The definition of the depth proposed in the aforementioned paper and has been modified for our model in this paper. In our case we have two wells and thus we will define two depths Eq. (\ref{depth}). The depth that is related to the difference between the potential energy of the upper index-1 saddle equilibrium
point (us) and the potential energy of the bottom well (bw) (see green line in panel A of  Fig. \ref{depth_flat}) and the depth that is related to the difference between the potential energy of the upper index-1 saddle equilibrium point (us) and the potential energy of the top well (tw) (see red line in  panel A of  Fig. \ref{depth_flat}). We notice that the depth evolves linearly, for both cases, as the asymmetry parameter takes values from 0 to 0.5.

\begin{equation}\label{depth}
    \text{depth} = V(x_{us},y_{us}) - V(x_{bw}, y_{bw}),\ (\text{depth} = V(x_{us},y_{us}) - V(x_{tw}, y_{tw}))
\end{equation}

Moreover in panel B of Fig. \ref{depth_flat} we show the difference between these two depths.

The flatness of our PES following  \cite{naikflat} is defined as

\begin{equation}\label{flatnes}
    \text{flatness} = \overline{||\frac{\partial V(x,y)}{\partial x}, \frac{\partial V(x,y)}{\partial y}||},\ (x,y) \in \Omega
\end{equation}

Thus the flatness is  the average of the Euclidean-norm of
the gradient of the potential energy function that is evaluated at discrete points in a bounded domain $\Omega$. In panel C of Fig. \ref{depth_flat} we show the flatness of our potential in two different bounded areas. The curve in green color represents the flatness calculated in the bounded area that the bottom well belongs, whereas the curve in red color represents the flatness calculated in the bounded area that the top well belongs. The flatness in both areas follows a quadratic law.

\subsection{Branching ratio and area of the lobes}

In this subsection we describe the main results of this paper. But first we describe the set up for our computations. The total energy of our system is fixed at $H = 0.1$. It is crucial to highlight that all the one thousand equally spaced initial conditions that we chose are lying in the line $x = - 0.005$. Moreover another condition for the trajectories is that $p_{x}>0$. In this way we consider only the trajectories that react (a trajectory reacts when it crosses the region
of the higher energy saddle (upper index-1 saddle) and approaches the lower energy saddle (lower index-1 saddle)). For the sake of completeness we note that we consider that a trajectory enters the top well if it crosses the line $y = 0.5$ and it enters the bottom well if it crosses the line $y = -0.5$. In both cases we stop the integration when the trajectories cross either of these lines.

\begin{figure}[htbp]
	\begin{center}
		A)\includegraphics[scale=0.165]{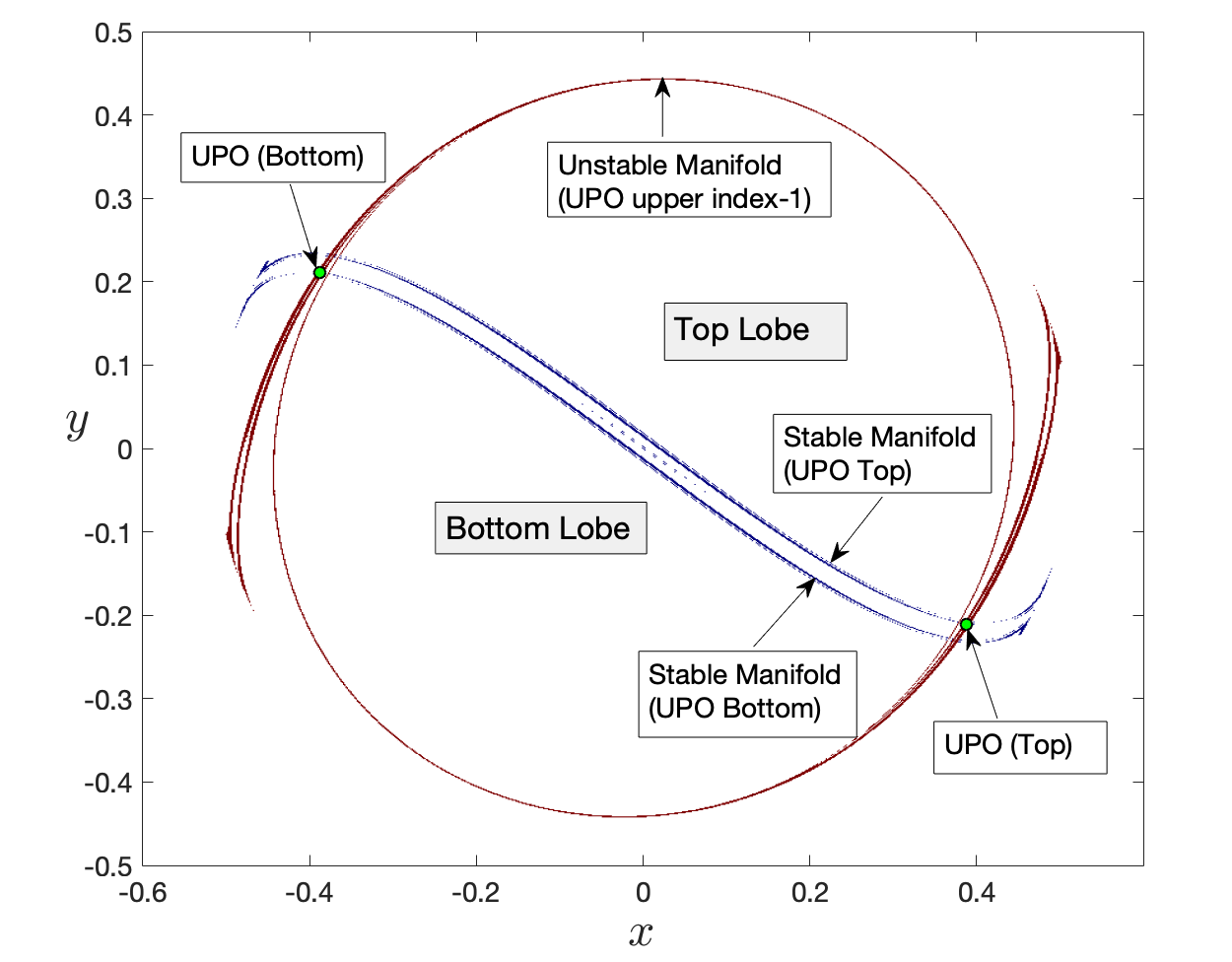} 
		B)\includegraphics[scale=0.20]{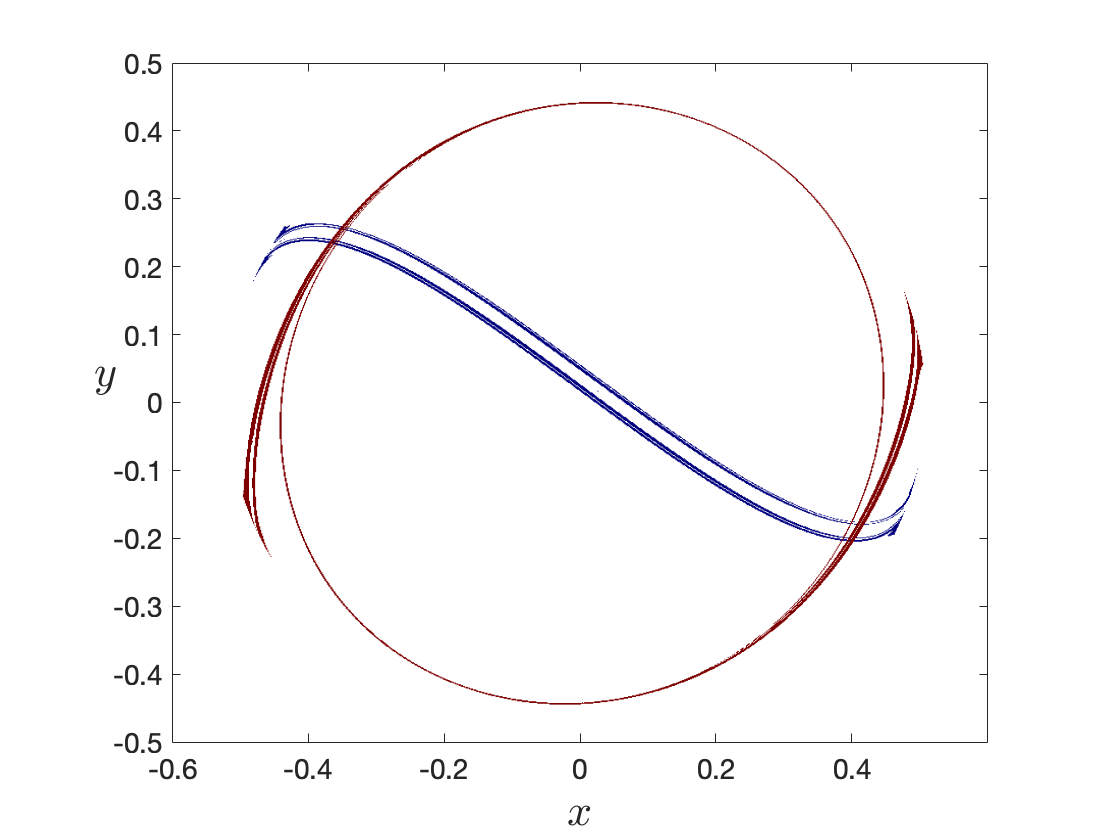}
		C)\includegraphics[scale=0.20]{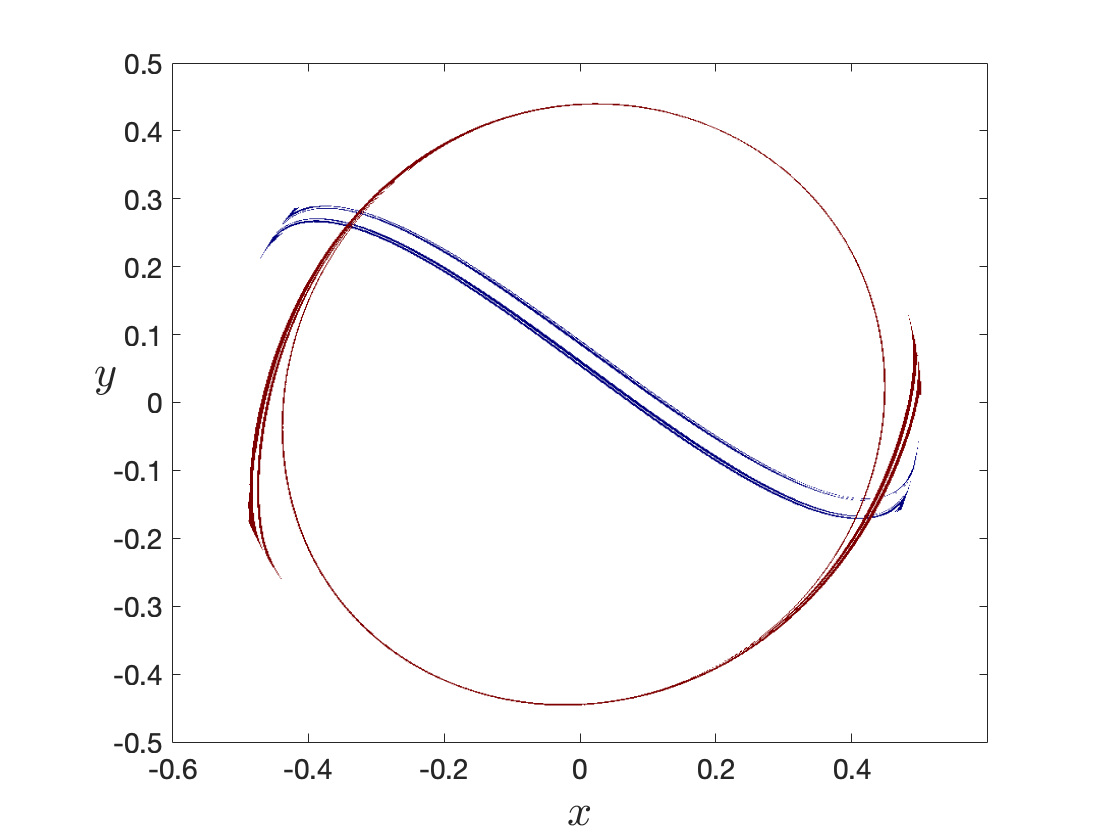}
	    D)\includegraphics[scale=0.20]{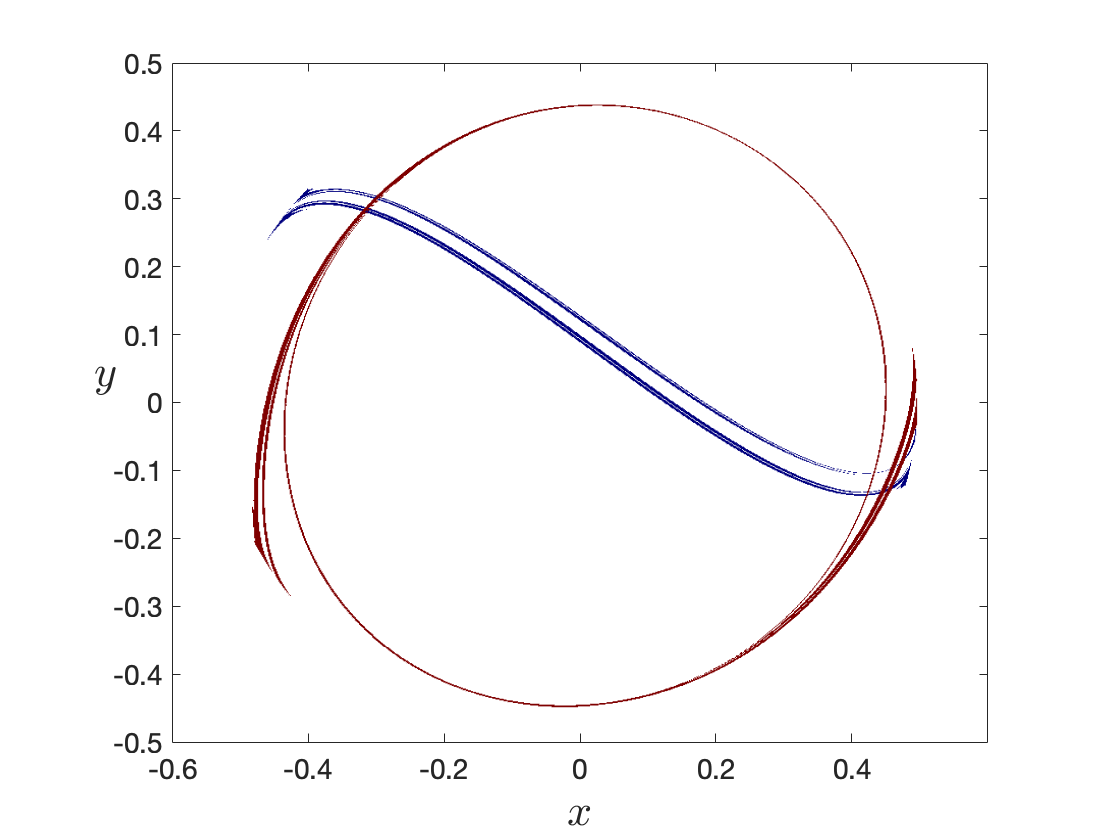}
	    E)\includegraphics[scale=0.20]{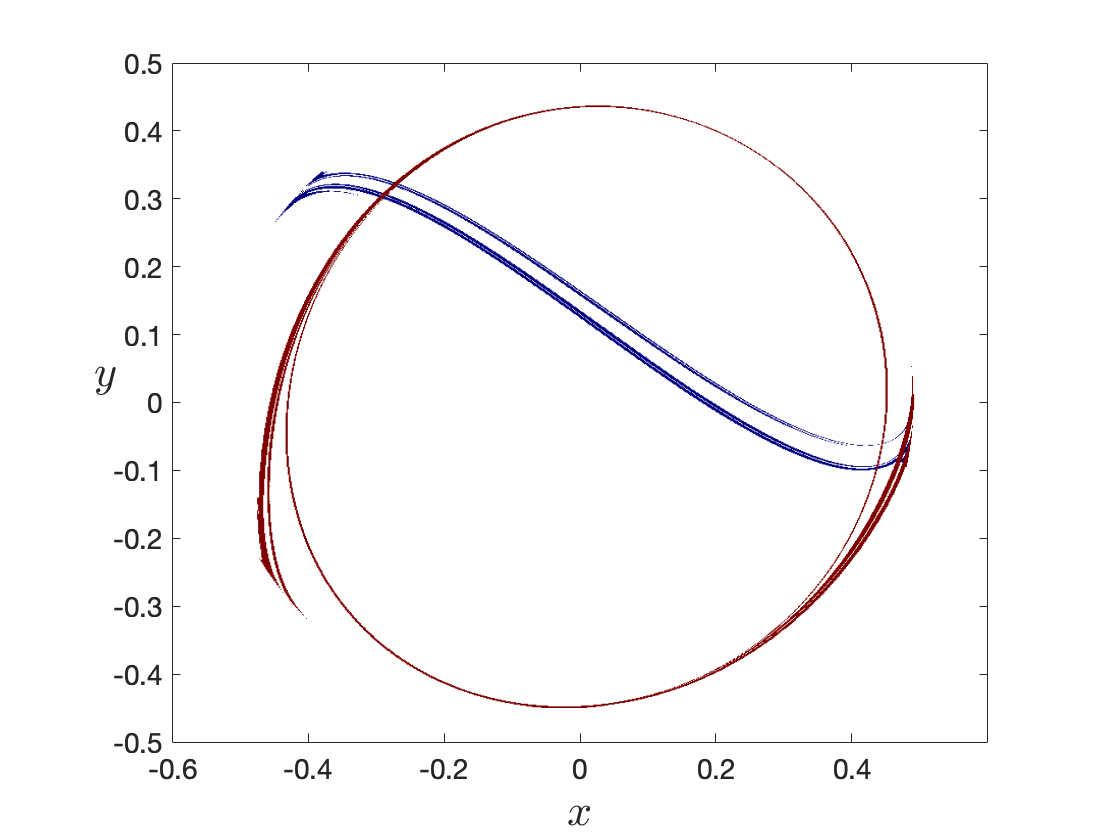}
	    F)\includegraphics[scale=0.20]{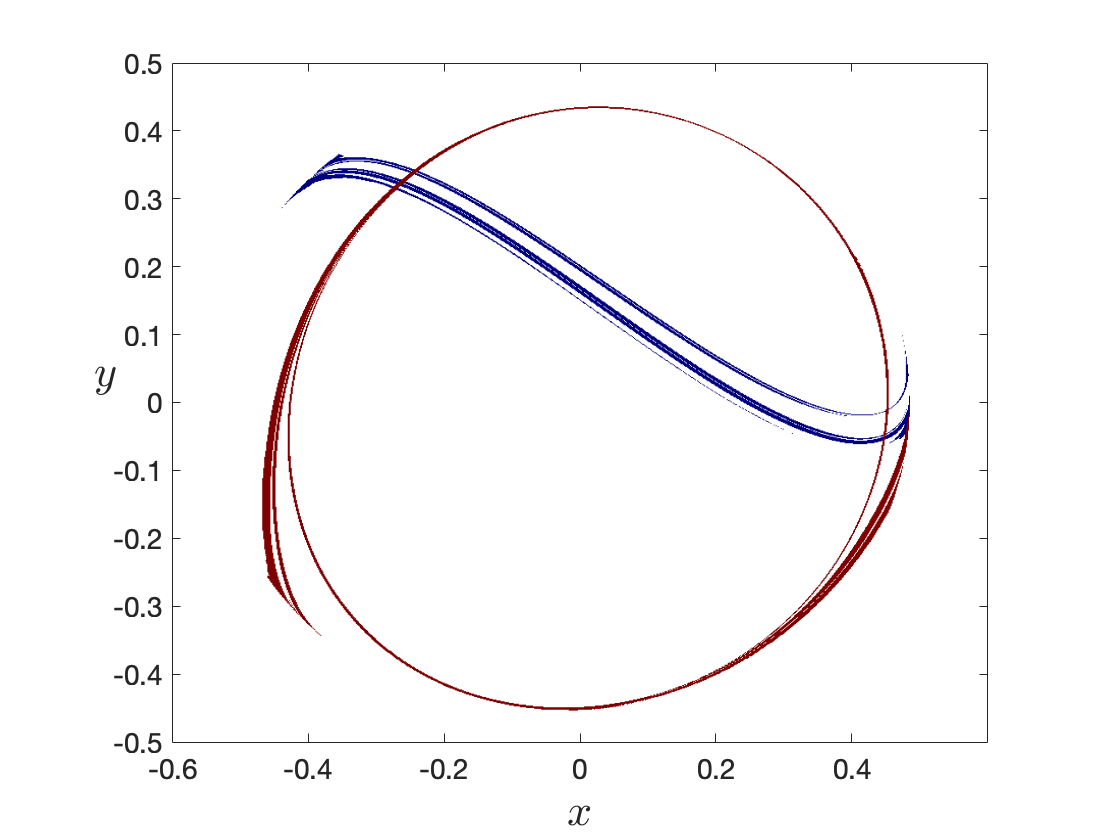}
	\end{center}
	\caption{Stable (blue) and unstable (red) invariant manifolds extracted from the gradient of the LD scalar field for different values of the asymmetry parameter $c$: A) $c = 0$ (symmetric case), B) $c = 0.1$, C) $c = 0.2$, D) $c = 0.3$, E) $c = 0.4$ and F) $c = 0.5$. }
	\label{fig:manifolds}
\end{figure}

In Fig. \ref{fig:manifolds} we present the stable (blue) and the unstable (red) manifolds that have been extracted from the gradient of the LD scalar fields, of the different unstable periodic orbits (UPO) in the system (for more details see paper \cite{Katsanikas2020}). Our focus for this paper are the two lobes that are associated with a heteroclinic intersection between the unstable manifold of the UPO of the upper index-1 saddle with the stable manifold of the top UPO (top lobe) and a heteroclinic intersection between the unstable manifold of the UPO of the upper index-1 saddle with the stable manifold of the bottom UPO (bottom lobe) respectively. As we described in \cite{Katsanikas2020} the top lobe  and the bottom lobe are responsible for the transport of the trajectories from the region of the upper index-1 saddle to the top and bottom wells, respectively. The manifolds have been calculated for six different values of the asymmetry parameter $c$. In panel A of Fig. \ref{fig:manifolds} we show the symmetric case where we can observe that the area of the lobes described above is identical for both lobes. When the value of the parameter $c$ is increasing the area of the lobes is no longer equal  and we see that the area of the top lobe decreases as the area of the bottom lobe increases. As the area of the lobe that is formed from the unstable invariant manifold of the unstable periodic orbits of the upper index-1 saddle with the stable invariant manifold of the top unstable periodic orbit decreases, the number of the trajectories that are trapped inside this lobe decreases too. Consequently, the ratio of the trajectories that visit the top well decreases. On the contrary, as the area of the lobe that is formed from the unstable invariant manifold of the unstable periodic orbits of the upper index-1 saddle with the stable invariant manifold of the bottom unstable periodic orbit increases, the number of the trajectories that are trapped inside this lobe increases too. Consequently, the ratio of the trajectories that visit the bottom well increases.

\begin{figure}[htbp]
	\begin{center}
		A)\includegraphics[scale=0.23]{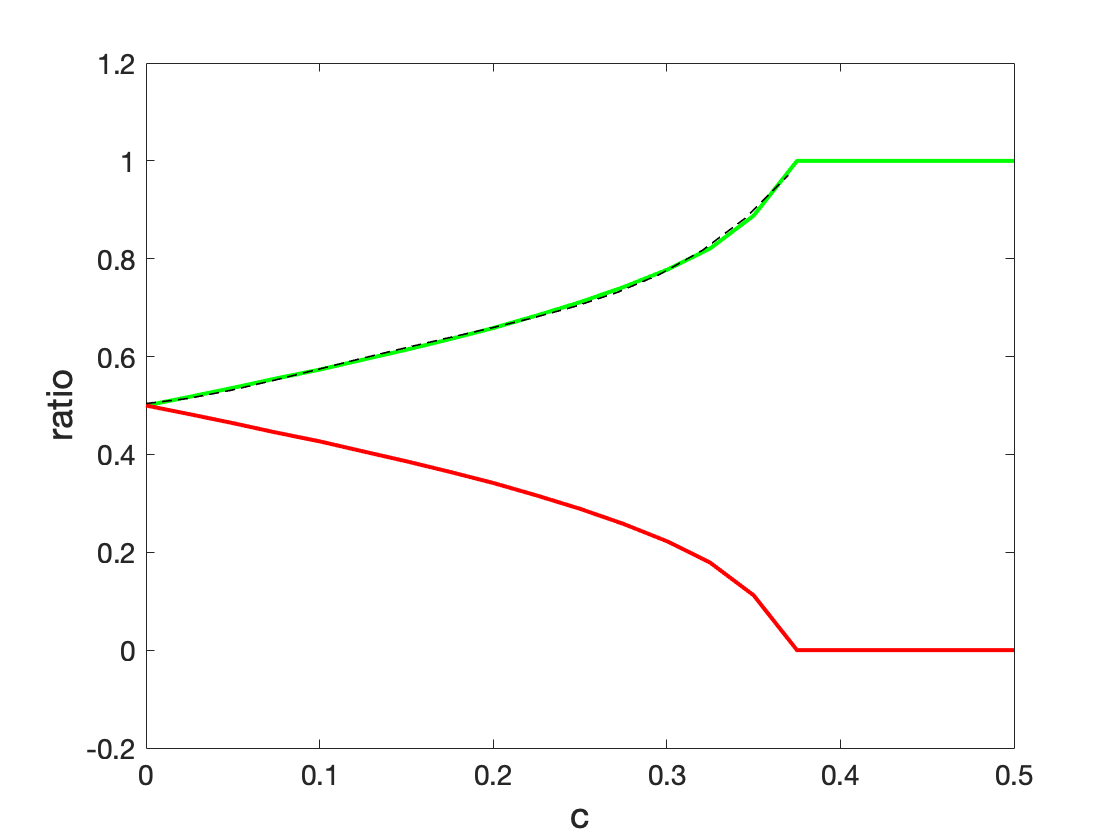} 
		B)\includegraphics[scale=0.23]{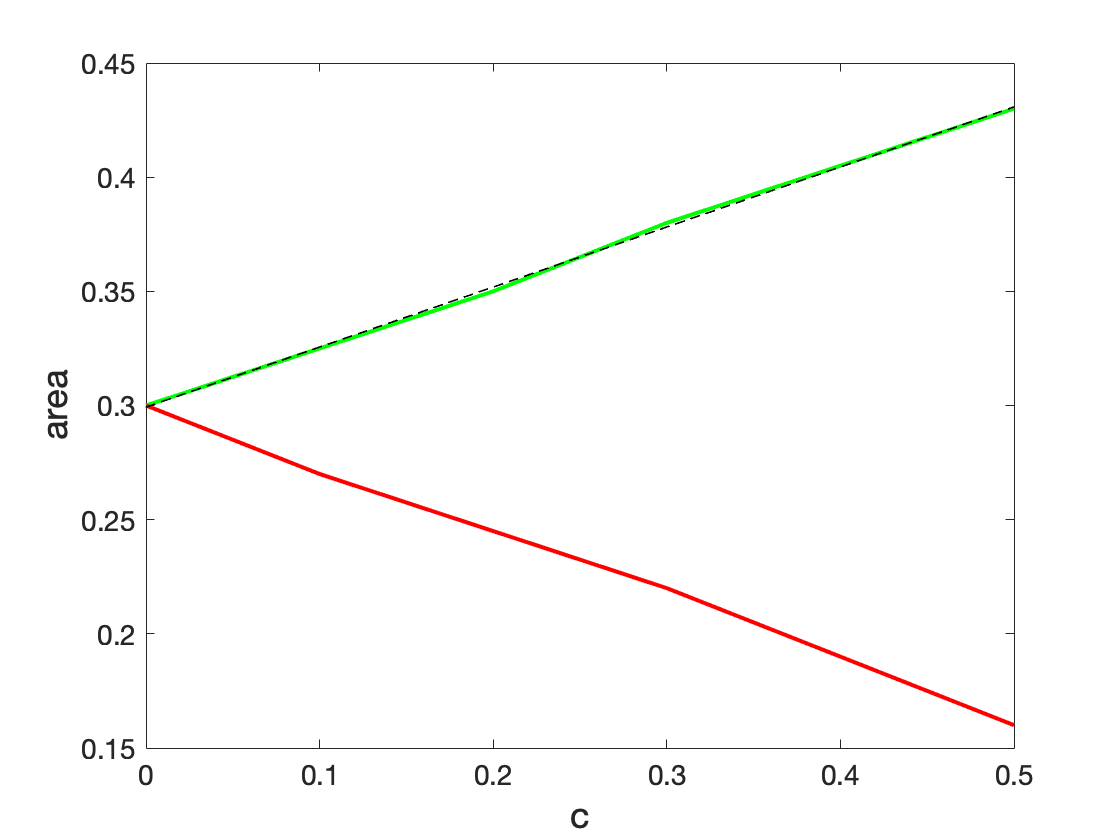}\\
		C)\includegraphics[scale=0.40]{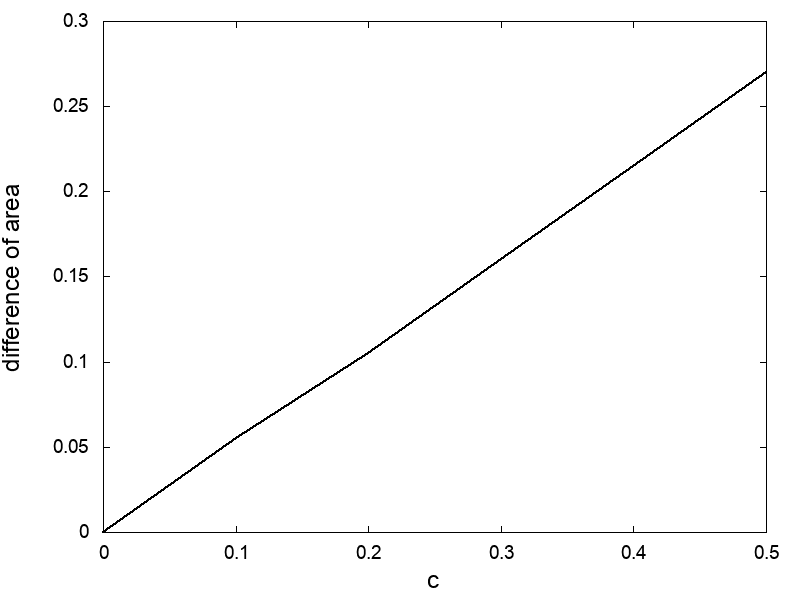}
	\end{center}
	\caption{A) The branching ratio of the trajectories that choose to visit the top or the bottom well as the parameter $c$ changes. On  top of the figure we present the approximation $r =  a_{1}c^{4} + a_{2}c^{3} +a_{3}c^{2}+a_{4}c +a_{5}$ where $a_{1} = 72.96, a_{2} = -40.89, a_{3} = 7.854,a_{4} = 0.2581, a_{5} = 0.504$ (black dashed line) for the ratio of the trajectories that enter into the bottom well.  B) The area of the lobes of the top and the bottom well as the parameter $c$ changes. On top of the figure we present the linear approximation $m = b_{1}c + b_{2}$ where $b_{1} =0.2629 , b_{2} = 0.2993$ for the bottom lobe area (black dashed line). C) The difference of the area of the lobes of the top well from the area of the lobe of the bottom well.}
	\label{fig:ratio}
\end{figure}

In the panel A of Fig. \ref{fig:ratio} we show the branching ratio of the trajectories that enter the top (red) or the bottom (green) well as the asymmetry parameter $c$ takes values between 0 and 0.5. We see in the Fig. \ref{fig:ratio} that the ratio of the trajectories that visit  the  top well is smaller than 
the ratio of the trajectories that visit the bottom well. This happens because, as we explained above, the top lobe  that is associated with the transport of the trajectories  from the region of the upper index-1 saddle to the region of the top well  becomes smaller than the bottom lobe  that is associated with the transport of the trajectories  from the region of the upper index-1 saddle to the region of the bottom well. Furthermore, we see that the two curves in the panel A of Fig. \ref{fig:ratio}  follow a quartic law (until a critical value of $c = 0.375$). The one curve that corresponds to the ratio of the trajectories that visit the top well  decreases, as the parameter $c$ increases, until it reaches zero (at $c=0.375$, nothing goes to the top well). On the contrary, the other curve  that corresponds to the ratio of the trajectories that visit the bottom well increases, as the parameter $c$ increases, until it reaches one (for $c=0.375$, all trajectories enter into the region of the bottom well). In panel B of Fig. \ref{fig:ratio}  we notice how the area of the top and bottom lobes evolves linearly as the asymmetry parameter changes in the range of values mentioned above and in a similar way with the evolution of the depth in panel A of Fig. \ref{depth_flat}. In  panel C of Fig. \ref{fig:ratio} we show that the difference of the areas of the top and bottom lobes  is linear. We observe that the  nonlinear evolution (quartic law) of the branching ratio is a result of the linear evolution  of the areas of the lobes that are responsible for the transport of the trajectories to the region of the wells. If our system was linear, this correlation would be linear (this means that a linear evolution of the areas should have as a result a linear evolution of the branching ratio). But our system is nonlinear and  the  correlation of the evolution of the areas with the  branching ratio is nonlinear and it follows a quartic law. Actually, this correlation is a measure of the nonlinearity of our system and how fast   our system will converge to  a state  in which all the trajectories visit only one well (in our case it converges for values of $c$ above  $0.375$).

\section{Conclusions}\label{Conclusions}
In this paper we studied an asymmetric PES which includes two sequential index-1 saddles, with one saddle having higher energy than the other at all values of the parameter $c$, that controls the asymmetry, and two wells, that have the same energy when $c$ is 0 and unequal energies as the parameter increases. In particular, as the parameter $c$ increases the energy of the top well increases while the energy of the bottom well decreases. Our goal is to show how the depth and the flatness of our PES is changing by varying the parameter $c$ and, moreover, to show the correlation between the area of the lobes (top and bottom) and the branching ratio of the trajectories that enter the top and bottom well while we vary the parameter $c$. 

For our analysis we use the method of Lagrangian Descriptors and trajectory calculations based on the set up that is mentioned in Section \ref{Results}.

Our main conclusions are:

\begin{enumerate}

    \item The depth as it is defined in Eq. (\ref{depth}) versus the parameter $c$  of the potential (that controls the asymmetry of our potential) follows a linear evolution while the flatness as it is defined in Eq. (\ref{flatnes}) versus the parameter $c$ in our potential follows a quadratic law.
    
    \item The evolution of the  branching ratio of the trajectories that  enter into the top and bottom well versus $c$ ( the parameter $c$ that  controls the asymmetry of our potential) is nonlinear and it follows a quartic law. This means that the ratio of the trajectories that enter to the one well increases fast (following a quartic law) until it becomes one. On the contrary, the ratio  of the trajectories that enter to the other well decreases fast (following a quartic law) until it becomes zero.
    
    \item The  area of the lobes, between the unstable invariant manifold of the unstable periodic orbits of the upper index-1 saddle and the stable invariant manifold of  the top and bottom unstable periodic orbits, versus the parameter $c$  of the potential evolves linearly and in a similar way to the evolution of the depth(s).  The difference of the areas between the two lobes  is also  linear.
    
    \item The linear  evolution   of the areas of  the two lobes that are responsible for the transport of the trajectories to the top or bottom well causes a quartic (nonlinear) evolution of the branching ratio of the trajectories that enter into the regions of the two wells. The correlation between the evolution of the areas and the evolution of the  branching ratio gives us  a measure of the nonlinearity of our system and  how fast our system will converge  to  a state  in which all the trajectories visit only one well. In our case this happens following a quartic law  and it converges at this state for values of the parameter $c$ above $0.375$.  
    
\end{enumerate}

\section*{Acknowledgements}
The authors acknowledge the support of EPSRC Grant No. EP/P021123/1 and MA also
acknowledges the support from the grant CEX2019-000904-S and IJC2019-040168-I funded by: MCIN/AEI/10.13039/501100011033.

\bibliography{SNreac}
\end{document}